\documentclass[12pt]{article}

\RequirePackage{a4}
\usepackage{epsfig}
\usepackage[utf8x]{inputenc}
\usepackage{verbatim}
\usepackage{lineno}
\usepackage{amsmath} 
\usepackage{amssymb}
\usepackage{hyperref}
\usepackage{url}
\usepackage{breakurl}

\newcommand{\bm}[1]{\mbox{\boldmath$#1$}}

\begin{document}

\title{The Gauss' Bayes Factor
}
\author{Giulio D'Agostini \\
Universit\`a ``La Sapienza'' and INFN, Roma, Italia \\
{\small (giulio.dagostini@roma1.infn.it,
 \url{http://www.roma1.infn.it/~dagos})}
}

\date{}

\maketitle

\begin{abstract}
In {\em Theoria motus corporum coelestium in sectionibus 
  conicis solem ambientum}
Gauss presents, as a theorem and with emphasis,
the rule to update the ratio of probabilities
of complementary hypotheses, in the light 
of an observed event which could be due to either of them.
Although he focused on {\em a priori} equally probable hypotheses, in order
to solve the problem on which he was interested in, the theorem can be
easily extended to the general case. But, curiously,
I have not been able to find references to his result in the literature.
\end{abstract}
\mbox{}
\begin{flushright}
  {\footnotesize \sl ``I play with a gentleman whom I do not know.\\
    He has dealt ten times, \\
    and he has turned the king up six times.\\ 
    What is the chance that he is a sharper?\\
    This is a problem in the  probability of causes.\\ 
    It may be said that it is the essential problem\\
    of the experimental method.''}\\
  (H. Poincaré)
\end{flushright}  

\section{Introduction}
As it is becoming rather well known, the only sound
way to solve what Poincaré called
{\sl ``the essential problem of the experimental method''}
is to tackle it using probability theory, as it should be rather
obvious for {\sl ``a problem in the  probability of causes''}.
The mathematical tool to perform what is also known
as `probability inversion' is called
{\em Bayes rule} (or theorem),
although due to Laplace, at least in one of the most common
formulations:\footnote{{\sl ``The greater the probability of an observed event
given any one of a number of causes to which that event may be attributed, 
the greater the likelihood of that cause $[$given that event$]$.
 The probability of the existence of any one of these causes
 $[$given the event$]$ is thus a fraction
whose numerator is the probability of the event given the cause,
and whose denominator is the sum of similar probabilities,
summed over all causes.
If the various causes are not equally probable 
{\em a priory}, it is necessary, instead of the probability of the event
given each cause, to use the product of this probability
and the possibility
of the cause itself.''}\cite{Laplace}
}
\begin{eqnarray}
       P(C_i\,|\,E,I) &=& \frac{P(E\,|\,C_i,I)\cdot P(C_i\,|\,I)}
       {\sum_k P(E\,|\,C_k,I)\cdot P(C_k\,|\,I)}\,,
                 \label{eq:BayesRule1}
\end{eqnarray}
where $E$ is the observed event and 
$C_i$ are its possible causes, forming a complete class
(i.e. exhaustive and mutually exclusive).
`$I$' stands for the background state of information,
on which all probability evaluations {\em do} depend
(`$I$' is often implicit, as it will be later in this paper,
but it is important to remember
of its existence).

Considering also an alternative cause $C_j$, 
the ratio of the two {\em posterior probabilities},
that is how the two hypotheses are re-ranked
in {\em degree of belief},
in the light of the observation $E$, is given by 
 \begin{eqnarray}
        \frac{P(C_i\,|\,E,I)}{P(C_j\,|\,E,I)} &=& 
         \frac{P(E\,|\,C_i,I)}{P(E\,|\,C_j,I)} \times 
         \frac{P(C_i\,|\,I)}{P(C_j\,|\,I)}\,,
        \label{eq:BayesRule2}
       \end{eqnarray}
 in which we have factorized the r.h. side into
 the {\em initial} ratio of probabilities
of the two causes (second term) and the updating factor
\begin{eqnarray}
             \frac{P(E\,|\,C_i,I)}{P(E\,|\,C_j,I)}\,,
             \label{eq:BF}
\end{eqnarray}
known as  {\em Bayes factor}, or  
 `likelihood ratio'.\footnote{The
alternative name {\em likelihood ratio} is  preferred
in some communities of researchers 
because numerator
and denominator of Eq.\,(\ref{eq:BF}) are called {\em likelihood}
by statisticians.
My preference to `Bayes factor' (or even {\em Bayes-Turing Factor}
\cite{Good-BTF,Banks-Good,GdA_GW}) is due to the fact
that, since in the common parlance `likelihood' and `probability'
are in practice equivalent, `likelihood ratio'
tends to generate confusion as it were the ratio of the probabilities
of the hypotheses of interest (and the value that maximizes the
`likelihood function' tends to be considered {\em by itself}
the most probable value).}
The advantage of
Eq.\,(\ref{eq:BayesRule2}) with respect to Eq.\,(\ref{eq:BayesRule1})
is that it highlights the two contributions to the {\em posterior}
ratio of the hypothesis of interest: 
the {\em prior}
probabilities of the `hypotheses',
on which there could be a large variety of opinions;
the ratio of the probabilities of the observed event,
under the assumption to each hypothesis of interest,
which can {\em often} be rather {\em intersubjective},
in the sense that there is usually a larger, or unanimous
consensus, if the conditions under they have been evaluated
(`$I$') 
are clearly stated and shared (and in critical cases
we have just to rely on the well argued and documented opinion
of experts.\footnote{For example the
  European Network of Forensic Science Institutes
 strongly recommends\,\cite{ENFSI} forensic scientists to report
 the `likelihood ratio' of the findings in the light
 of the hypothesis of the prosecutor and the hypothesis
 of the defense, abstaining to assess which hypothesis
 they consider more probable, task left to the judicial
 system (but then I have strong worries, shared by other
 researchers, about the ability of the members of
 judicial system of making
 the proper use of such a quantitative information!).
 To those interested on the details of how
 this {\em Guideline} can be turned into practice,
 a {\em Coursera} offered by the University of Lausanne
 is recommended\,\cite{CourseraLausanne}
})

Recently, going after years through  
the third section  of the second
`book' of Gauss' {\em Theoria motus corporum coelestium in sectionibus 
  conicis solem ambientum}\,\cite{Gauss, Gauss_trasl},
of which I had read with the due care only the part
in which the {\em Prince Mathematicorum} derives in his
peculiar way what is presently known as the Gaussian (or `normal')
error function,
I have realized that Gauss had also illustrated, a few pages before,
a theorem on how to update the probability ratio of two
alternative hypotheses, based on experimental observations.
Indeed the theorem is not exactly Eq.(\ref{eq:BayesRule2}), because
it is only formulated for the case in which
$P(C_i\,|\,I)$ and $P(C_j\,|\,I)$ are equal,
but the reasoning Gauss had setup would have led naturally
to the general case. It seems that he focused
into the sub-case of {\em a priori} equally likely hypotheses
just because he had to apply his result to a problem
in which he consider the values to be inferred
{\em a priori} equally likely 
({\sl ``valorum harum incognitarum
  ante illa observationes aeque probabilia fuisse''}).

But let us proceed in order.

\section{Probability of observations vs probability of the values
of physical quantities}
The third section of `book 2' of the Gauss' tome\,\cite{Gauss} 
is dedicated
to {\sl ``the determination of an orbit satisfying as nearly as possible
any number of observations whatever''}.\footnote{All English quotes 
are taken from the C.H. Davis translation\,\cite{Gauss_trasl}.}
After `articles'\,\footnote{This publication is divided into two `books',
each of them subdivided in four `sections'. Then the entire
text is divided
in {\em numeri} (translated into `articles' by Davies\,\cite{Gauss_trasl})
running through the `books'. In particular, Section 3 of Book 2,
consisting of 22 printed pages in the original Latin edition,
contains `articles' 172 to 189. 
} 172-174, which introduce
the specific problem of evaluating the
elements of an orbit from the measurements
of geocentric quantities related to those elements, 
with article 175 Gauss {\em ascends}\,\footnote{``\ldots {\em ad disquisitionem
    generassimam in omni calculi ad philosophiam naturalem applicatione
    fecundissima ascendemus.}''}
to methodological
issues of general interest for the Sciences:
\begin{quote}
  {\sl ``let us leave our special problem,
    and enter upon a very general discussion and
    one of the most fruitful in every application
    of the calculus to the natural philosophy.''}
\end{quote}
The general problem is how to determine the $\mu$ unknown quantities
$p$, $q$, $r$, $s$, etc. (e.g. the elements of the orbit
of a planet or a comet) and evaluate the functions $V_i$ of these variables
from $\nu$ measurements $V_{m_i}$ (e.g. the geocentric quantities
of that celestial body
measured at different times):\footnote{For
  the reader' convenience (hopefully)
  the functions are called here, except when they appear in quotes,
  $V_1$, $V_2$, $V_3$,  etc.,
  and the measured values $V_{m_1}$, $V_{m_2}$, $V_{m_3}$, etc.,
  while Gauss uses  $V$, $V'$,  $V''$\ldots and
  $M$, $M'$,  $M''$ \ldots, respectively.}
\begin{eqnarray}
 V_i(p,q,r,s,\dots)  & \xrightarrow [\mbox{\it measured as\ }]{} & V_{m_i}\,, \label{eq:M-V}
\end{eqnarray}
or, indicating the set of unknown quantities by $\bm{\theta}$,
that is $\bm{\theta} = \{p,q,r,s,\dots\}$, we can rewrite
Eq.\,(\ref{eq:M-V}) as 
\begin{eqnarray}
  V_i(\bm{\theta})  & \xrightarrow [\mbox{\it measured as\ }]{}& V_{m_i} \,. \label{eq:M-V-theta}
\end{eqnarray}
The most interesting case, Gauss explains, is when
$\nu > \mu$. Being over-determined, 
this case has a solution only if $V_{m_i}$ are affected by experimental
errors, described by the probability density
function\footnote{Note how Gauss
  simply speaks of `probabilities',
  obviously meaning {\em probability
    density functions}, as clear from the
  use he makes of them: {\sl ``the probability to be assigned
  to each error $\Delta$ will be expressed by a function of $\Delta$
  that we shall denote as $\varphi\Delta$''} -- a few lines later it is
  clear that Gauss had in mind a 'pdf' since, when he wrote
 {\sl ``the probability generally, that the error lies between $D$ and $D'$,
  will be given by the integral $\int\!\varphi\Delta\,\mbox{d}\Delta$ extended
  from $\Delta= D$ to  $\Delta= D'$''. 
  $[\,$Note also that in the case a function had only one argument, 
    parentheses were not used. Therefore  $\varphi\Delta$ stands
    for  $\varphi(\Delta).]$}\label{fn:pdf}}
(pdf)
$\varphi$, that  in article 177 will came out
to be the well known Gaussian function.\footnote{For
  an account of Gauss' derivation in modern notation
  see Sec. 6.12 of Ref.~\cite{BR}
  (some intermediate steps needed to reach the solution
are sketched in
\url{http://www.roma1.infn.it/~dagos/history/Gauss_Gaussian.pdf}).}
Therefore,\footnote{As clarified in footnote \ref{fn:pdf}, it is
  clear that in the following quotes the generic term
  {\sl ``probability''} stands for {\em probability density function}.
}
\begin{quote}
  {\sl
    ``Supposing, therefore, any determinate system of the values of the
    quantities $p$, $q$, $r$, $s$, etc.,
    the probability that the observation would give for $V$ the value $M$
    will be expressed by $\varphi(M\!-\!V)$, substituting in $V$ for
    $p$, $q$, $r$, $s$, etc., their values;  in the same manner
    $\varphi(M'\!-\!V')$,   $\varphi(M''\!-\!V'')$, etc, will express the
    probability that observation would give the values $M'$, $M''$, etc.
    of the functions  $V'$, $V''$, etc. Wherefore, since we are authorized
    to regard all observations as event independent of each other, the
    product
    $$\mbox{} \hspace{2.0cm}\varphi(M\!-\!V)\,\varphi(M'\!-\!V')\,
    \varphi(M''\!-\!V'')\,\mbox{etc,}\, = \,\Omega
    \hspace{2.0cm}\mbox{\rm (G1)}$$
    will express the expectation or probability that all those values
    will result together from observation.''
  }  
  \end{quote}
What Gauss calls $\Omega$ is thus the {\em joint pdf}
of the  differences $V_{m_i}-V_i$
given a precise set of values for the
physical quantities of interest, which we would rewrite as
\begin{eqnarray}
  f(\bm{V_m}\!-\!\bm{V}\,|\,\bm{\theta})
  &=& \prod_i\varphi(V\!_{m_i}\!-\!V_i\,|\,\bm{\theta})
  \label{eq:Omega}
\end{eqnarray}
where $\bm{V_m}$ and  $\bm{V}$ stand for the set of observations
and of functions.\footnote{The
  reason why in the argument appears the differences and not simply
  the observed values is that for Gauss $\varphi()$ was the error
  function, i.e. `probability density function' of the errors
  (see also footnote \ref{fn:pdf}). Since, later in `article' 177,
  the function  $\varphi()$ will become the `Gaussian' error function,
  we could rewrite directly the joint pdf of the observations
  in modern notation
  as\\ \mbox{}\hspace{3.0cm}
  $f(\bm{V_m}\,|\,\bm{\theta}) = \prod_i\frac{1}{\sqrt{2\pi}\,\sigma}
  \exp\left[-\frac{\left(V_{m_i}-V_i(\bm{\theta})\right)^2}{2\sigma^2}\right]$.
  \label{fn:gaussiana}
} 

Article 175 ends so with the expression of the
joint probability of the observations given any set
of values of the quantities of interest, that is a problem
in {\em direct probabilities}:
$$\bm{\theta} \xrightarrow [\mbox{\it deterministic link\  }]{} \bm{V}
\xrightarrow [\mbox{\it probabilistic link\  }]{}\bm{V_m} $$
Article 176 begins with what we could call nowadays
a `Bayesian manifesto':
\begin{quote}
  {\sl
    ``Now in the same manner as, when any determinate values whatever
    of the unknown quantities being taken, a determinate probability corresponds,
    previous to observations, to any system of values of the
    functions $V$, $V'$, $V''$, etc;
    so, inversely, after determinate values of the functions have resulted
    from observation, a determinate probability will belong
    to every system of values of the
    unknown quantities, from which the values of the
    functions could possibly have resulted.''    
  }  
\end{quote}
That is, in our notation, 
as when we assume {\sl ``determinate values''}
of the  physical quantities 
we are interested in the joint pdf of the
values that will be observed,
\begin{eqnarray}
  f(\bm{V_m}\!-\!\bm{V}\,|\,\bm{\theta})\,,  \label{eq:f_M_given_theta}
\end{eqnarray}
similarly, once the observations have been made,
we are interested in the joint pdf of the values
of the physical  quantities,
\begin{eqnarray}
  f(\bm{\theta}\,|\,\bm{V_m}\!-\!\bm{V})\,. \label{eq:f_theta_given_M}
\end{eqnarray}
The question is now how to go from Eq.\,(\ref{eq:f_M_given_theta})
to Eq.\,(\ref{eq:f_theta_given_M}), reasoning ``inversely''.

\section{Updating the probabilities of hypotheses}
\begin{figure}
\centering{\fbox{\epsfig{file=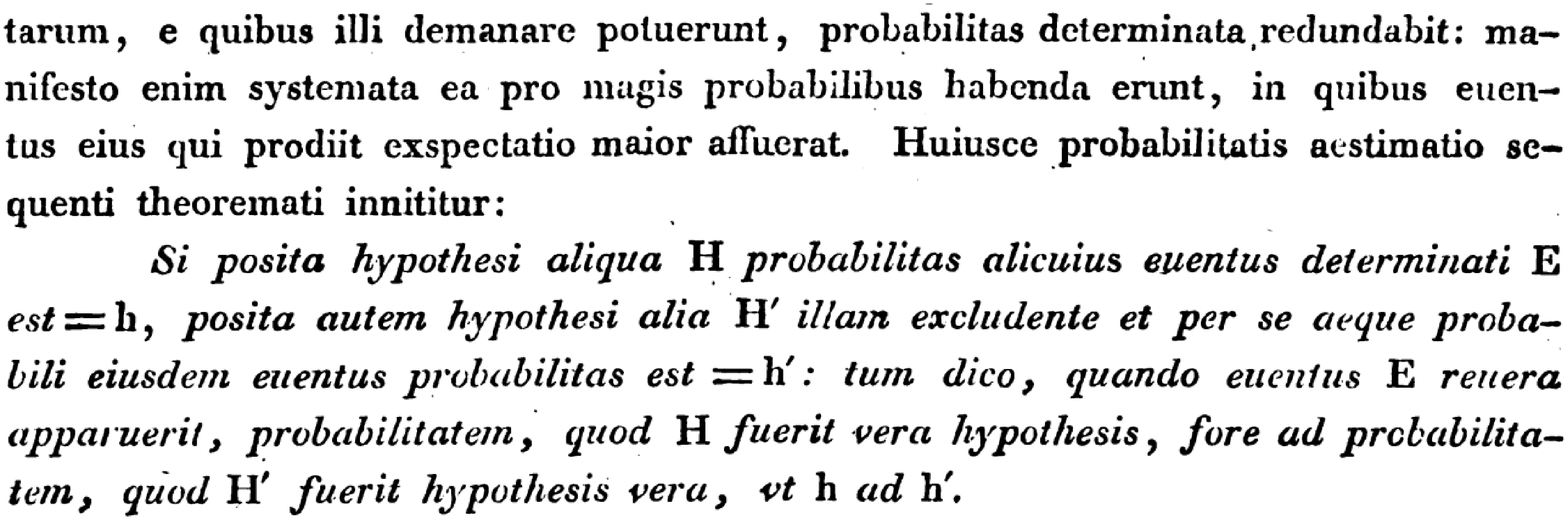,clip=,width=\linewidth}}}
\caption{\small \sf Extract of {\em Theoria motus corporum\ldots}\cite{Gauss}
  in which Gauss enunciates {\em his} theorem on how to
  update probability ratios of incompatible hypotheses in the
  light of an experimental observation. Note ``{\bf tum dico}'' (``than I say'').
}
\label{fig:GaussBF}
\end{figure}

We are finally at the core of the problem. Let Gauss speak:
\begin{quote}
  {\sl
    ``For, evidently, those systems will be regarded as the more probable
    in which the greater expectation had existed of the event
    which actually occurred.
    The estimation of this probability rests upon the following theorem:\\
    \mbox{}\ \ \ \ {\it If, any hypothesis  H being made, the probability of any
      determinate event E is h, and if, another hypothesis H' being made excluding
      the former and equally probable in itself, the probability of the same
      event is h': then I say, when the event E has actually occurred,
      that the probability that H was the true hypothesis, is to the
      probability that H' was the true hypothesis, as h to h'.}''
  } (Italic original, also put in evidence in the text as a quote
  -- see Fig.~\ref{fig:GaussBF}.)  
\end{quote}
In modern notation:
\begin{eqnarray}
  P(E\,|\,H) &=& h \nonumber \\
  P(E\,|\,H') &=& h' \nonumber  \\
  \frac{P(H\,|\,E)}{P(H'\,|\,E)} &=& \frac{P(E\,|\,H)}{P(E\,|\,H')}\,,
  \hspace{0.5cm}\mbox{{\bf  if} $P_0(H) = P_0(H')$}\,. \label{eq:GaussBF}
\end{eqnarray}
There are no doubts that Gauss presents this result as
original ({\sl ``then I say''}, in Latin {\em tum dico}), although
it might be curious that it did not refer to results by Laplace,
who had been writing on {\em probabilities of causes} more than
thirty years before\footnote{The historian of statistics
  Stephen Stigler refers to Laplace's 1774 {\em Mémoire}
  as {\sl ``arguably the most influential article this field}
  [{\em  mathematical statistics}$^{(*)}$] {\sl to appear before 1800,
  being the first widely read presentation of inverse probability
  and its application to both binomial and location
  parameter estimation.''}\,\cite{Stigler_Laplace}
  (note that in this reference there is no mention to Gauss).\\
  $^{(*)}$\,As far as I know, neither Gauss nor Laplace were
  using the word `statistics', but they were talking about probability.
}\,\cite{Laplace_PC}.
(For comparison, a few pages later, in article 177, Gauss
acknowledges Laplace
for having calculated the integral needed to normalize the `Gaussian'
distribution.) It is also curious the fact that Gauss starts saying
that {\sl ``evidently, those systems will be regarded as the more probable
    in which the greater expectation had existed of the event
    which actually occurred''}, considering
    thus {\sl ``evident''}
    what is presently known as `maximum likelihood principle',
    but then taking care of proving it as a theorem (under
    the well stated assumption of initially equally probable hypotheses).

The reasoning upon which the theorem is proved
is based on an inventory of equiprobable cases. This might seems
to limit the  application to
situations in which this inventory is in practice feasible, like
in games of cards and of dice. Instead, this was the 
way of reasoning of those times
to partition the space of possibilities, as it is clear from
the use that Gauss makes of his result, certainly not limited
to simple games.
\begin{figure}
\centering{\fbox{\epsfig{file=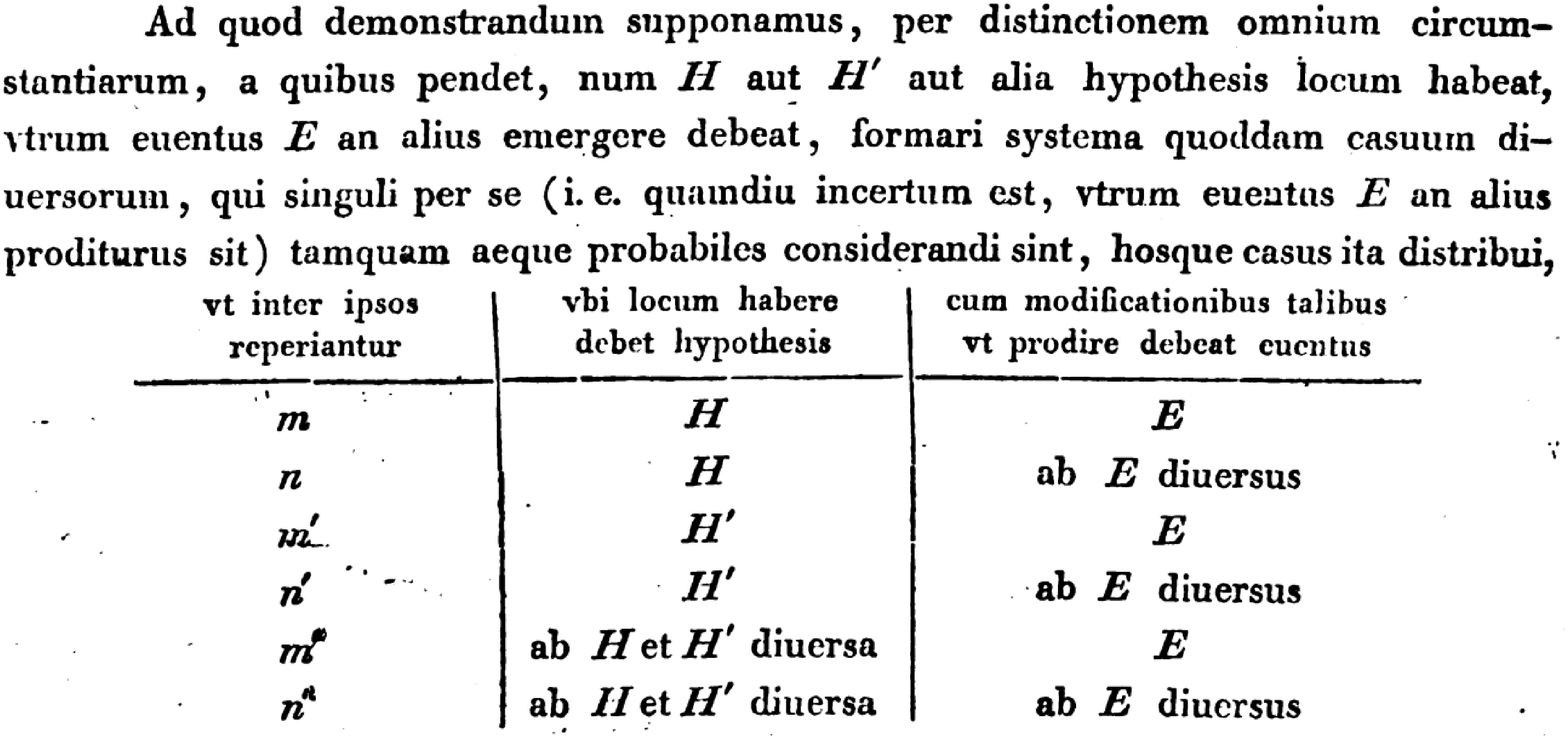,clip=,width=\linewidth}}}
\caption{\small \sf Partition of the space of possibilities
  as it appears in the original work of Gauss\,\cite{Gauss}.
  The English translations of the three columns are\,\cite{Gauss_trasl}:
  ``that among them may be found'';
  ``in which should be assumed the hypothesis''; ``in such a mode
  as would give occasion to the event''. Then:
  ``ab $E$ diuersus''\,$=$\,``different from $E$'';
  ``ab $H'$ et $H'$ diuersa''\,$=$\,``different from $H$ and $H'$\,''. 
}
\label{fig:H:h_table}
\end{figure}
Figure \ref{fig:H:h_table} shows the original version of such a partition.
The six numbers of the first column, normalized to their sum, provide
the following probabilities:
\begin{eqnarray}
  P(E\cap H) &=& \frac{m}{m+n+m'+n'+m''+n''} \nonumber\\
  && \nonumber \\
  P(\overline{E}\cap H) &=& \frac{n}{m+n+m'+n'+m''+n''}\nonumber \\
  && \nonumber \\ 
  P(E\cap H') &=& \frac{m'}{m+n+m'+n'+m''+n''} \nonumber\\
  && \nonumber \\
  P(\overline{E}\cap H') &=& \frac{n'}{m+n+m'+n'+m''+n''} \nonumber\\
  && \nonumber \\
  P(E\cap \overline{H\cup H'}) &=& \frac{m''}{m+n+m'+n'+m''+n''} \nonumber\\
    && \nonumber \\
  P(\overline{E}\cap \overline{H\cup H'}) &=& \frac{n''}{m+n+m'+n'+m''+n''}\nonumber 
\end{eqnarray}
The probabilities which enter the proof are those
of the $H$ and $H'$
\begin{eqnarray}
  P(H) &=& \frac{m+n}{m+n+m'+n'+m''+n''}\label{eq:P0H}\\
  && \nonumber \\
  P(H') &=& \frac{m'+n'}{m+n+m'+n'+m''+n''} \label{eq:P0H'}
\end{eqnarray}
and those of the event $E$ given
either hypothesis:
\begin{eqnarray}
  P(E\,|\,H) &=& \frac{m}{m+n} = h \label{eq:P_E_H}\\
   && \nonumber \\
  P(E\,|\,H') &=& \frac{m'}{m'+n'} = h' \label{eq:P_E_H'}
\end{eqnarray}
The probability of $H$ is modified by the observation of $E$
observing that, with
reference to Eqs.\,(\ref{eq:P0H}) and (\ref{eq:P0H'}),
\begin{quote}
{\sl ``after the event is known, when the cases $n$, $n'$, $n''$}
 {\sl  disappear
from the number of possible cases, the probabilities of the
same hypothesis will be
$$ \frac{m}{m+m'+m''}\,; $$
in the same way the probability of the
hypothesis $H'$ before and after the event,
respectively, will be expressed by
$$ \frac{m'+n'}{m+n+m'+n'+m''+n''}\ \ \ \mbox{and}\ \ \ 
\frac{m'}{m+m'+m''}\,: $$
since, therefore, the same probability is assumed for the hypotheses
$H$ and $H'$ before the event is known, we shall have
$$\mbox{} \hspace{4.5cm}
m+n\  =\  m'+n' \,, \hspace{4.5cm}\mbox{\rm (G2)}$$
hence the truth of the theorem is readily inferred.''}
\end{quote} 
That is, in our notation,
\begin{eqnarray}
  P(H\,|\,E) &=& \frac{m}{m+m'+m''}\nonumber \\
   && \nonumber \\
  P(H'\,|\,E) &=& \frac{m'}{m+m'+m''}\,,\nonumber
\end{eqnarray}
from which
\begin{eqnarray}
  \frac{P(H\,|\,E)}{P(H'\,|\,E)} &=& \frac{m}{m'}.\nonumber
\end{eqnarray}
Using then Eqs.\,(\ref{eq:P_E_H}) and (\ref{eq:P_E_H'}),
yielding $m=(m+n)\cdot P(E\,|\,H)$
and \\  $m'=(m'+n')\cdot P(E\,|\,H')$, we obtain
\begin{eqnarray}  
  \frac{P(H\,|\,E)}{P(H'\,|\,E)}
  &=& \frac{P(E\,|\,H)\cdot (m+n) } {P(E\,|\,H')\cdot (m'+n') }
  \label{eq:preBF}
\end{eqnarray}
Applying finally the condition (G2), theorem (\ref{eq:GaussBF})
is proved.

In reality, it is easy to see that,
being
$$\frac{m+n}{m'+n'} =\frac{P(H)}{P(H')}\,,$$
Eq.\,(\ref{eq:preBF})
contains the most general case
\begin{eqnarray}
  \frac{P(H\,|\,E)}{P(H'\,|\,E)} &=&
   \frac{P(E\,|\,H)} {P(E\,|\,H')}\cdot \frac{P(H)}{P(H')}\,. \nonumber
\end{eqnarray} 
But Gauss contented himself with the sub-case of initially
probable hypotheses.
Why? The reason is most likely that he focused on the
inference
of the unknown values of the physical quantities of interest,
that he assumed {\em a priori} equally likely, 
a very reasonable  assumption for this kind of inferences,
if we compare the prior knowledge
with the information provided by observations
(see e.g. Ref.\,\cite{BR}).

\section{Application to the inference of unknown values of physical quantities}
In fact, immediately after the proof of {\em his} theorem, Gauss continues:
\begin{quote}
  {\sl
    ``Now, so far as we suppose that no other data exist for the determination
    of the unknown quantities besides the observations
    $V=M$, $V'=M'$, $V''=M''$ etc., and, therefore, that all
    systems of values of these unknown quantities were equally probable previous
    to the observations, the probabilities, evidently,
    of any determinate system subsequent to the observations
    will be proportional to $\Omega$. This is to be understood
    to mean that the probability that the values of the unknown quantities
    lie between the infinitely near limits $p$ and $p+\mbox{d}p$,
    $q$ and $q+\mbox{d}q$,  $r$ and $r+\mbox{d}r$,
    $s$ and $s+\mbox{d}s$,, etc. respectively, is expressed
    by
    $$\ \ \ \ \ 
    \lambda\,\Omega\,\mbox{d}p\,\mbox{d}q\,\mbox{d}r\,\mbox{d}s\,\cdots,\ 
    \mbox{etc.}, \hspace{1.5cm}\mbox{\rm (G3)}
    $$
    where the quantity $\lambda$ will be a constant quantity independent
    of $p$, $q$, $r$, $s$, etc.:
    and, indeed,, $1/\lambda$ will, evidently, be the value of the integral
    of order $\nu$ ,
     $$\int^\nu\!\Omega\,\mbox{d}p\,\mbox{d}q\,\mbox{d}r\,\mbox{d}s\,\cdots,\ 
    \mbox{etc.},  \hspace{1.5cm}\mbox{\rm (G4)}
    $$
    for each of the variables  $p$, $q$, $r$, $s$, etc, extended
    from the value $-\infty$ to the value  $+\infty$.''
    }
\end{quote}
As we can see, it is well stated the assumption of `flat priors', as we use
to say nowadays (with the original words of Gauss, in Latin:
{\em ``valorum harum incognitarum
  ante illa observationes aeque probabilia
  fuisse''}).\footnote{It is clear that what is unknown are the
  numeric values of the quantities
  and not the `quantities' themselves, at it could seem from the English
  translation, because
  in that case
  there would be little to infer.} 

It is, instead, less clear how he uses the result of his
theorem (the quote at the beginning of this section follows
immediately the end of the proof of the theorem,
with no single word in between). The implicit intermediate step
is
\begin{equation}
  P(H\,|\,E) \propto  P(E\,|\,H)\,,
\end{equation}
extended to set of continuous uncertain values 
(`uncertain vector') $\bm{\theta}$
as
\begin{equation}
  P(\bm{\theta}\,|\,\mbox{data}) \propto  P(\mbox{data}\,|\,\bm{\theta})\,.
\end{equation}
Then, remembering that $\Omega$ was  the joint pdf of
the observations $[$see Eq.\,(G1)$]$, which we have rewritten
in more  compact notation as Eq.\,(\ref{eq:Omega}),
we have
\begin{eqnarray*}
  f(\bm{\theta}\,|\,\bm{V_m}-\bm{V}) & \propto&
  f(\bm{V_m}-\bm{V}\,|\,\bm{\theta}) \\
  \mbox{or}\hspace{4.0cm} && \\
  f(\bm{\theta}\,|\,\bm{V_m}-\bm{V}) & = & \lambda\cdot
  f(\bm{V_m}-\bm{V}\,|\,\bm{\theta})\,, 
\end{eqnarray*}
being $\lambda$ just the normalization constant, i.e.\footnote{Note:
  the reason of using in the formulae $\bm{V_m}-\bm{V}$, instead
than just $\bm{V_m}$ is simply due to the way Gauss wrote the error function,
but, obviously, this function could be redefined and $\bm{V}$ would disappear
from the above equations, then getting for example
$ f(\bm{\theta}\,|\,\bm{V_m}) \propto f(\bm{V_m}\,|\,\bm{\theta}) $,
as we would write it nowadays
(see also footnote \ref{fn:gaussiana}).}
\begin{eqnarray}
 \frac{1}{\lambda}  =
  \int_{\mathbb{R}^\nu}f(\bm{V_m}-\bm{V}\,|\,\bm{\theta})\,\mbox{d}\bm{\theta}\,.
\end{eqnarray}

\section{Conclusions}
Reading Gauss' work, there are no doubts that the
{\em Prince Mathematicorum}  had clear ideas on how to tackle
inverse probability problems, i.e.  
what goes presently under the name {\em Bayesian inference}.
In particular, he presented as original 
what is now called Bayes factor, i.e the factor to update
the {\em odds} in favor of an hypothesis with respect to
the alternative one, in the light of a new observation.
However, it is curious that, as far as I could find,
this result is not acknowledged
in the current literature.
For example, his name appears only once in the Sharon Mcgrayne
rather comprehensive book on the history of Bayesian
reasoning\cite{McGrayne}, as being cited by Enrico Fermi,
who was teaching his students data analysis methods
derived from {\em his} Bayes'
theorem.\footnote{Indeed Ref.\,\cite{McGrayne} cites
Ref.\,\cite{Fermi_Bayes}, writing which I had realized
that Gauss was using a `Bayesian reasoning',
but I had at that time completely skipped
the `details' in which he derived, as a theorem,
the rule to update the ratio of probabilities of hypotheses,
subject of this paper.}

At this point a long discussion could follow on the question
if Gauss could be classified as a {\em Bayesian} and why,
later on in his book, he did not proceed applying 
consistently the probabilistic reasoning he had setup, getting
the joint probability distribution of the values of the
orbital elements given the observed geocentric measurements, but
he derived, instead, the {\em least square} method to get
(relatively) simple formulae for the {\em most probable values}
(this aim was clearly stated).
And all this in the same text, just a few pages after, and not
in a later stage of his life.

Well, I am not an historian,
and therefore I can only state my impressions based
on a limited amount of reading. Gauss appears 
in the section of the book upon which this modest note
is based not only as the genius he is famous to be, but also
a very practical scientist going straight to his goals. 
Trying to set a multi-dimensional inference to write
down the joint pdf of parameters of a non-linear problem
and exploiting it at best,
something that we can do nowadays,
thanks to unprecedented computing power and novel
mathematical methods, would have just been a waste of time
two centuries ago. We have also seen that he didn't even
care to state the general rule to update probability ratios, which
would have required just a couple of lines of text,
because he had in mind a problem for which the priors were
reasonable `flat'. Moreover, he was also well aware
of the practical meaning and limits of the mathematical functions,
as when, later in the same section, he commented in `article' 177 
on the {\sl ``defect''} of {\em his error function},
because {\sl ``the function just found cannot, it is true, express
rigorously the probabilities of the errors''}.
Indeed, the `error function' $\varphi()$ was not specified
up to the end of 'article' 176.
Only in the following article he showed
that a good candidate for it was, under well stated conditions,
\ldots \  the Gaussian, a function having the  {\sl ``defect''}
of contemplating values ranging from minus infinity to plus infinity.
Then  other interesting articles follow,\footnote{For example
  he derived the formula of the weighted average, stressing its
  importance to use it, instead of the individual
  values\,\cite{GdA_curious}.}
but I don't
want to spoil you the pleasure of the
reading.\footnote{Historical French and German translations
  are also available\,\cite{Gauss_F,Gauss_D}.}

Finally, someone might be intrigued about what Gauss meant
by {\em probability}. {\sl ``Probabilitas''}.
What else?\footnote{Although the formal theory
of probability has only been developed in the last few centuries,
the noun {\em probabilitas}, the adjective
{\em probabilis} and especially its comparative {\em probabilior}
(`more probable'), playing a fundamental
role in probabilistic reasoning,
were used in Latin with essentially the
same meaning we assign to them in ordinary language.
For example, a recent `grep' through the Cicero texts collected
in The Latin Library\,\cite{LatinLibrary} resulted in 
105 words containing `probabil'. And it is rather popular
the Cicero's quote {\sl ``Probability is the very guide
  of life''}\,\cite{CiceroProbability}, although such a sentence
does not appears verbatim in his texts, but it is a digest of his
thought (see e.g. {\em De Natura Deorum}, Liber Primus,
nr. 12 \cite{LatinLibrary}).
}


\end{document}